\input amstex

\NoBlackBoxes

\documentstyle{amsppt}

\topmatter

\title{On representability of algebraic functions by radicals}\endtitle

\author Askold Khovanskii\endauthor

\affil{Department of Mathematics, University of Toronto, Toronto, Canada}\endaffil


\thanks{This work was partially supported by the Canadian Grant No. 156833-17.}\endthanks

\abstract{This preprint is dedicated to a self contained simple proof of the classical criteria for  representability of  algebraic functions of several complex variables by  radicals. It also contains a  criteria for  representability of  algebroidal functions by composition of single-valued analytic functions and radicals, and a  result related to the 13-th Hilbert problem. This preprint is an extended version of the author's 1971 paper. It is written as a part of the comments to a new edition (in preparation) of the classical book ``Integration in finite terms'' by J.F.~Ritt.  }\endabstract

\keywords{algebraic function,  solvability by radicals, 13-th Hilbert problem  }\endkeywords


\endtopmatter
\document

Consider an algebraic equation
$$ P_n y^n+P_{n-1}y^{n-1}+\dots +P_0=0, \tag 1 $$ whose coefficients  $P_n,\dots, P_0$ are polynomials  of $N$ complex variables $x_1,\dots, x_N$.
Camille Jordan discovered that  the  Galois group of the equation (1) over the field $\Cal R$ of rational functions of $x_1,\dots,x_N$ has a topological meaning (see theorem 3 below): it is isomorphic to the {\it monodromy group} of the equation (1).

According to   Galois theory, the equation (1) is solvable
by  radicals over the field $R$ if and only if
its Galois group is solvable. If the equation (1) is irreducible it defines a multivalued algebraic function $y(x)$.
Galois theory and Theorem 3 imply  the  following criteria for representability of an algebraic function by radicals, which consists of  two statements:

1) {\it  If the  monodromy group of an algebraic function  $y(x)$ is
solvable, then $y(x)$ is representable by radicals.}

2) {\it \it  If the  monodromy group of an algebraic function  $y(x)$ is not
solvable, then $y(x)$ is not representable by radicals.}

 We reduce the first statement  to linear algebra (see Theorem 10 below) following the book [4].

We prove the second statement  topologically without using Galois theory. Vla\-di\-mir Igorevich Arnold found the first topological proof of this statement [1]. We use another topological approach (see Theorem 15 below) based on the paper [3].  This paper contains  the first result of  topological Galois theory [4] and it gave a hint for its further development.

\subhead{1. Monodromy group and Galois group}
\endsubhead
Consider the  equation (1). Let $\Sigma \subset \Bbb C^N$ be the hypersurface defined by  equation $P_n J=0$, where $P_n$ is the leading coefficient and $J$ is the discriminant of the equation (1). The  {\it  monodromy group } of the equation (1)  is the group of all  permutations of  its solutions  which are induced by
motions  around the singular set $\Sigma $ of the equation (1). Below we discuss this definition more precisely.

At a point $x_0\in \Bbb C^N\setminus \Sigma $ the set $Y_{x_0}$   of all germs   of analytic  functions  satisfying the equation (1) contains exactly $n$ elements, i.e. $Y_{x_0}=\{y_1,\dots,y_n\}$. Indeed, if  $P_n(x_0)\neq 0$ then for    $x=x_0$  the equation (1) has $n$ roots counted with multiplicities. If in addition   $J(x_0)\neq 0$ then all these roots are simple. By the implicit function theorem each simple root can be extended to a germ of a regular function satisfying the equation  (1).

Consider a closed curve  $\gamma $ in
$\Bbb C^N\setminus \Sigma $ beginning and ending at the point $x_0$. Given a germ $y\in Y_{x_0}$ we can continue it along the loop $\gamma$ to obtain another germ $y_\gamma\in Y_{x_0}$.
Thus each such loop $\gamma$   corresponds to a permutation $S_{\gamma}:Y_{x_0}\rightarrow Y_{x_0}$
 of the set $Y_{x_0}$ that maps a germ $y\in Y_{x_0}$ to the germ $y_{\gamma}\in Y_{x_0}$. It is  easy to see that the map $\gamma\rightarrow S_\gamma$ defines a homomorphism   from the fundamental
group $\pi _1(\Bbb C^N\setminus \Sigma,x_0)$ of the domain $\Bbb C^N\setminus \Sigma $ with the base point $x_0$ to the group $S(Y_{x_0})$ of permutations. The {\it
monodromy group} of the equation (1) is the image of
the fundamental group   in the group
$S(Y_{x_0})$ under this homomorphism.

\remark{Remark} Instead of the point $x_0$ one can choose any other point $x_1\in \Bbb C^N\setminus \Sigma$. Such a choice will not change the    monodromy group  up to an isomorphism. To fix this isomorphism one can choose any curve $\gamma:I\rightarrow \Bbb C^N\setminus \Sigma$ where $I$ is the segment $ 0\leq t\leq 1$ and $\gamma(0)=x_0$, $\gamma(1)=x_1$ and identify each germ $y_{x_0}$ of solution of (1) with its continuation $y_{x_1}$ along $\gamma$.

Instead of the hypersurface $\Sigma$ one can choose any bigger algebraic hypersurface $D$, $\Sigma\subset D\subset \Bbb C^N$. Such a choice will not change the    monodromy group: one can slightly  move a curve $\gamma \in \pi_1(\Bbb C^N\setminus \Sigma, x_0)$ without changing the map $S_\gamma$ in such a way that $\gamma$ will not intersect $D$.
\endremark

The field of rational functions of $x_1,\dots, x_N$ is isomorphic to the field $\Cal R$ of germs of rational functions at the point $x_0\in \Bbb C^N\setminus \Sigma $.
Consider the field extension  $\Cal R \{y_1,\dots,y_n\}$ of $\Cal R$ by the  germs  $ y_1,\dots y_n$ at $x_0$  satisfying the equation (1).

\proclaim {Lemma 1} Every permutation $S_\gamma$ from the monodromy group can be uniquely extended to an automorphism of the   field $\Cal R\{ y_1,\dots
,y_n\}$ over the field $\Cal R$.
\endproclaim

\demo {Proof} Every element $f\in \Cal R \{y_1,\dots,y_n\}$ is a rational function of $x, y_1,\dots,y_n$. It can be continued meromorphically along  the curve
$\gamma\in \pi_1(\Bbb C^m\setminus \Sigma,x_0)$  together with  $y_1,\dots,y_n$.
This continuation gives the required
automorphism, because the continuation preserves the arithmetical
operations  and every rational function returns back
to its original  values (since it is a single-valued valued function). The automorphism is unique because the extension is generated by $y_1,\dots, y_n$.
\enddemo

By definition the {\it  Galois  group} of the equation (1)  is the group of all automorphisms of  the   field $\Cal R\{ y_1,\dots,y_n\}$ over the field $\Cal R$. According to Lemma 1 the monodromy group  of the equation (1) can be considered as a subgroup of its Galois group. Recall that by definition a multivalued function $y(x)$ is  {\it algebraic} if all its meramorphic germs satisfy the same algebraic equation over the field of rational functions.

\proclaim {Theorem 2} A germ $f\in \Cal R\{y_1,\dots,y_n\}$ is fixed under the monodromy action if and only if  $f\in \Cal R$.
\endproclaim
\demo {Proof} A germ $f\in \Cal R\{y_1,\dots,y_n\}$ is fixed under the monodromy action  if and only if $f$   is a germ of a single  valued function. The field $\Cal R\{y_1,\dots,y_n\}$ contains only germs of algebraic functions.  Any single valued algebraic function is a rational function.
\enddemo

According to  Galois theory Theorem 2 can be formulated in the following way.

\proclaim {Theorem 3} The monodromy group of  the equation (1) is isomorphic to the Galois group of the equation (1) over the field $\Cal R$.
\endproclaim

Below we will not rely on Galois theory. Instead  we will use Theorem 2 directly.

\proclaim{Lemma 4} The monodromy group acts on the set $Y_{x_0}$ transitively if and only if the equation (1) is irreducible over the field of rational functions.
\endproclaim

\demo{Proof} Assume that there is a proper subset $\{y_1, y_2,\dots y_k\}$   of $Y_{x_0}$  invariant under the monodromy action. Then the basic symmetric functions $r_1=y_1+\dots +y_k$, $r_2=\sum_{i<j}y_i y_j,$ $\dots$, $r_k=y_1\cdot\dots\cdot y_k$ belong to the field $\Cal R$. Thus $y_1, y_2,\dots y_k$ are solutions of the degree $k<n$ equation $y^k-r_1y^{k-1}t+\dots (-1)^kr_k=0$. So the equation (1) is reducible. On the other hand if the equation (1) can be represented as a product of two equations over $\Cal R$ then their roots belong to  two complementary  subsets of $Y_{x_0}$ which are  invariant under the monodromy action.
\enddemo

\proclaim{Corollary 5} An irreducible  equation (1)  defines a multivalued algebraic function $y(x)$ whose set of germs at $x_0\in \Bbb C^N\setminus \Sigma $ is the set $Y_{x_0}$ and whose monodromy group coincides with the monodromy group of the equation (1).
\endproclaim

Theorem 3, Corollary 5 and Galois theory immediately imply the following result.

\proclaim {Theorem 6} An algebraic function whose monodromy group is solvable can be represent by rational functions using the arithmetic operations and radicals.
\endproclaim
A stronger version of Theorem 6  can be proven  using linear algebra (see Theorem 10 in the next section).

\subhead { 2. Action of  solvable groups and representability
by radicals}\endsubhead
In this section, we prove that if a finite solvable group $G$ acts on a $\Bbb C$-algebra $V$ by
automorphisms, then  all elements of $V$ can be expressed by the elements of the invariant subalgebra $V_0$ of $G$ by taking radicals and adding.

 This construction of a representation by radicals is based on linear algebra. More precisely we use the following well known statement:  any finite abelian
group of linear transformations of a finite-dimensional vector space over  $\Bbb C$
can be diagonalized in a suitable basis.

 \proclaim{Lemma 7}  Let $G$ be a finite abelian group of order $n$ acting by automorphisms
on  $\Bbb C$-algebra $V$.
Then every element of the algebra $V$ is representable as a sum of  elements
$x_i \in V $,  such that $x^n_i$
lies in the invariant subalgebra $V_0$ of $G$, i.e., in
the fixed-point set of the group $G$.
\endproclaim

\demo {Proof} Consider a finite-dimensional vector subspace $L$ in the algebra $V$ spanned by
the $G$-orbit of an element $x$. The space $L$ splits into a direct sum $L = L_1+\dots+L_k$
of eigenspaces for all operators from $G$. Therefore, the vector $x$ can
be represented in the form $x = x_1+\dots + x_k$, where $x_1,\dots, x_k$ are eigenvectors
for all the operators from the group. The corresponding eigenvalues are $n$-th roots of
unity. Therefore, the elements $x^n_1,\dots, x^n_k$ belong to the invariant subalgebra $V_0$.
\enddemo

\definition {Definition} We say that an element $x$ of  algebra $V$ is an $n$-th root of an
element~$a$ if $x^n=a$.
\enddefinition

We can now restate Lemma 7 as follows: every element $x $ of the algebra $V$
is representable as a sum of $n$-th roots of some elements of the invariant subalgebra.

\proclaim {Theorem 8} Let $G$ be a finite solvable group of order $n$ acting by automorphisms
on  $\Bbb C$-algebra $V$. Then
every element $x$ of the algebra $V$ can be obtained from the elements of the invariant
subalgebra $V_0$ by takings $n$-th roots  and summing.
\endproclaim

We first prove the following simple statement about an action of a group on a set.
Suppose that a group $G$ acts on a set $X$, let $H$ be a normal subgroup of $G$, and denote by
$X_0$  the subset of $ X$ consisting of all points fixed under the action~of~$G$.

\proclaim {Lemma 9} The subset $X_H$ of the set $X$ consisting of the fixed points under the action of the normal subgroup  $H$ is invariant under the action of $G$. There is a
natural action of the quotient group $G/H$ on the set $X_H$ with the fixed-point set~$X_0$.
\endproclaim

\demo {Proof} Suppose that $g\in G$, $h \in H$. Then the element $g^{-1}hg$ belongs to the normal subgroup $H.$ Let $x\in X_H$. Then $g^{-1}hg (x)= x$, or $h(g(x))= g(x)$, which means
that the element $g(x)\in X $ is fixed under the action of the normal subgroup $H$. Thus
the set $X_H$ is invariant under the action of the group $G$. Under the action of $G$ on
$X_H$, all elements of $H$ correspond to the identity transformation. Hence the action
of $G$ on $X_H$ reduces to an action of the quotient group $G/H$.
\enddemo

We now proceed with the proof of Theorem 8.
\demo {Proof (of Theorem 8)} Since the group $G$ is solvable, it has a chain of nested
subgroups $G = G_0\supset \dots\supset G_m = e$ in which the group $G_m$ consists of the
identity element $e$ only, and every group $G_i$ is a normal subgroup of the group
$G_{i-1}$. Moreover, the quotient group $G_{i-1}/G_i$ is abelian.
Let $V_0 \subset \dots\subset V_m =V$ denote the chain of invariant subalgebras of the
algebra $V$ with respect to the action of the groups $G_0,\dots,G_m$. By Lemma 9
the abelian group $G_{i-1}/G_i $ acts naturally on the invariant subalgebra $V_i $, leaving
the subalgebra $V_{i-1}$  pointwise fixed. The order $m_i$ of the quotient group $G_{i-1}/G_i$
divides the order of the group $G$. Therefore, Lemma 7 is applicable to this
action. We conclude that every element of the algebra $V_i$ can be expressed with the
help of summation and $n$-th root extraction by the elements of the algebra $V_{i-1}$.
Repeating the same argument, we will be able to express every element of the
algebra $V$ by the elements of the algebra $V_0 $ using a chain of summations and $n$-th
root extractions.
\enddemo

\proclaim {Theorem 10} An algebraic function whose  monodromy  is solvable can be represented by rational functions by root extractions and summations.
\endproclaim

\demo {Proof} One can prove Theorem 10 by applying Theorem 8 to the monodromy action by automorphisms on the extension $\Cal R\{y_1,\dots,y_n\}$ with the field of invariants $\Cal R$.
\enddemo

\subhead {3. Topological obstruction to  representation by radicals}\endsubhead
Let us introduce some notation.

By $G^m$ we denote the $m$-th commutator of the group $G$. For any  $m\geq $ the group $G^m$ is a normal subgroup in $G$.

By $F(D,x_0)$ we denote the fundamental group  of the domain $U=\Bbb C^N\setminus D$ with the base point $x_0\in U$, where $D$ is an algebraic hypersurface in $\Bbb C^N$.

Let $H(D,m)$ be the covering space  of the domain $\Bbb C^N\setminus D$ corresponding to the subgroup $F^m(D, x_0)$ of the fundamental group $F(D,x_0)$.

We will say that an algebraic function  is an {\it  $R$-function} if it  becomes a single-valued function on some covering $H(D,m)$.

\proclaim{Lemma 11} If $m_1\geq m_2$ and $D_1\supset D_2$ then there is a natural projection
$\rho:H(D_1,m_1)\rightarrow H(D_2,m_2)$. Thus if a function $y$  becomes a single-valued function on $H(D_2,m_2)$ then it certainly becomes a single-valued function on $H(D_1,m_1)$.
\endproclaim

\demo {Proof} Let $p_*: F(D_1,x_0)\rightarrow F(D_2,x_0)$ be the homomorphism induced by the embedding $p: \Bbb C^N\setminus D_1\rightarrow \Bbb C^N\setminus D_2$. Lemma 11 follows from  the following  obvious  inclusions: $p_*^{-1}[F^{m_2}(D_2,x_0)]\subset F^{m_2}(D_1,x_0)$ and $F^{m_2}(D_1,x_0)\subset F^{m_1}(D_1,x_0)$.
\enddemo

\proclaim{Lemma 12} If $y_1$ and $y_2$ are $R$-function then $y_1+y_2$, $y_1-y_2$, $y_1\cdot y_2$ and  $y_1/y_2$  also are $R$-functions.
\endproclaim

\demo {Proof} Assume that $R$-functions $y_1$ and $y_2$ become single-valued functions on the coverings $H(D_1,m_1)$ and  $H(D_2,m_2)$. By Lemma 11 the  functions $y_1$,$y_2$ become single-valued on the covering $H(D,m)$ where $D=D_1 \bigcup D_2$ and $m=\max (m_1, m_2)$. Thus the functions $y_1+y_2$, $y_1-y_2$, $y_1\cdot y_2$ and  $y_1/y_2$  also become single-valued on on the covering $H(D,m)$.
The proof is completed since  $y_1+y_2$, $y_1-y_2$, $y_1\cdot y_2$ and  $y_1/y_2$ are algebraic functions.
\enddemo

\proclaim {Lemma 13} Composition of an $R$-function with the degree $q$ radical  is an $R$-function.
\endproclaim

\demo {Proof} Assume that the function $y$ defined by (1) is  $R$-function which becomes a single-valued function on the covering $H(D_1,m)$.  Let $D_2\subset \Bbb C^N$ be the hypersurface, defined by the equation $P_n P_0=0$, where $P_n$ and $P_0$ are the leading coefficient and the constant term of the equation (1). According to Lemma 11 the function $y$ becomes a single-valued function on the covering $H(D,m)$ where $D=D_1\bigcup D_2$. Let $h_0\in H(D,m)$  be a point whose image under the natural projection $\rho:H(D,m)\rightarrow \Bbb C^N\setminus D$ is the point $x_0$. One can identify the fundamental groups  $\pi_1 (H (D,m),h_0)$ and $F^m(D, x_0)$.

By  definition of  $D_2$ the function $y$ never equals to zero or to infinity on $H(D,m)$. Hence $y$  defines a map $y:H(D,m)\rightarrow \Bbb C\setminus \{ 0\}$.    Let $y_*:\pi_1 (H (D,m),h_0)\rightarrow \pi_1(\Bbb C\setminus \{0\}, y(h_0))$ be the induced homomorphism of the fundamental groups. The group $\pi_1 (H (D,m),h_0)$ is identified with the group $F^m(D,x_0)$ and the group $\pi_1(\Bbb C\setminus \{0\}, y(h_0))$ is isomorphic to $\Bbb Z$. So $\ker y_*\subset F^{m+1}(D,x_0)$. Thus all loops from the group $y_*(F^{m+1}(D,x_0))$ do not wind around the origin $0\in \Bbb C$.  Hence any germ of  $y^{1/q}$ does not change its value after continuation along a loop from the group $F^{m+1}(D,x_0)$. So $y^{1/q}$ is a single- valued function on $H(D,m+1)$. The proof is completed since $y^{1/q}$ is an algebraic function.
\enddemo

\proclaim {Lemma 14} An algebraic function $y$ is an $R$-function if and only if its monodromy group  is solvable.
\endproclaim
\demo {Proof} Assume that $y$ is defined by (1). Let $D$ be the hypersurface $P_n J=0$ where $P_n$ is the leading coefficient and $J$ is the discriminant of (1). Let  $M$ be  the monodromy group of $y$. Consider   the natural homomorphism $p:F(D,x_0)\rightarrow M$. If $M$ is solvable then for some natural number $m$ the $m$-th commutator of $M$ is the identity element $e$. The function $y$ becomes  single-valued on the covering $H(D,m)$ since $F^m(D,x_0)\subset p^{-1}(M^m)=p^{-1}(e)$. Conversely,  if $y$ is a single-valued function on some covering $H(D,m)$ then $p(F^m(D,x_0))=e$. But $p(F^m(D,x_0))=M^m$. Thus the monodrogy group $M$ is solvable.
\enddemo

\proclaim {Theorem 15} If an algebraic function  has unsolvable monodromy group
then it can not by represented by compositions of rational functions and radicals
\endproclaim

\demo {Proof} Lemma 12 and Lemma 13 show that the class of $R$-functions is closed under arithmetic operations and compositions with radicals. Lemma 14 shows that the monodromy group of any $R$-function is solvable.
\enddemo

\subhead {4. Compositions of analytic functions and radials}\endsubhead
In this section   we describe a class of multivalued functions in a domain $U\subset \Bbb C^N$
 representable by composition of single-valued analytic functions and radicals.

 A multivalued function $y$ in $U$  is called an {\it algebroidal function} in $U$ if it satisfies an irreducible equation
$$
y^n + f_{n-1}y^{n-1} + \dots +f_0=0\tag 2
$$
whose coefficients $f_{n-1},\dots, f_0$ are analytic functions in $U$. An algebroidal function could be considered as a continuous multivalued function in $U$ which has finitely many values.
\proclaim {Theorem 16 ([2], [3])} A multivalued  function $y$ in the domain $U$   can be  represented by composition of radicals and single valued analytic functions  if and only  $y$ is an algebroidal  function in $U$ with solvable monodromy group.
\endproclaim

To prove the ``only if" part one can repeat the proof of Theorem 15 replacing coverings over domains $\Bbb C^N\setminus D$ by coverings over domains $U\setminus \tilde D$ where  $\tilde D$ is an analytic hypersurface in $U$.

  To prove Theorem 16 in the opposite direction one can use Theorem 8 in the same way as it was used in the proof of Theorem 10.

\subhead {5. Local representability}\endsubhead In this section   we describe a  a local version of Theorem 16.

Let $y$ be an algebroidal function in $U$ defined by (2).  One can localize the equation (2)  at any point $p\in U$, i.e. one can replaced the coefficients $f_i$ of the equation (2)  by their germs at $p$. After such a localization  the equation (2) can became reducible, i.e. it can became representable as a  product of irreducible equations. Thus an algebroidal functions  $y$ in  arbitrary small neighborhood of a point $p$ defines several  algebroidal functions, which we will call {\it ramified germs of $y$ at $p$.} For a ramified germ of $y$ at $p$ the monodromy group is defined (as the monodromy group of an algebroidal function in an arbitrary small  neighborhood of the point $p$).

A ramified germ of  an algebroidal function $y$   of one variable $x$ in a neighborhood of a  point $p\in \Bbb C^1$ has a simple structure: its monodromy group is a cyclic group $\Bbb Z/m  \Bbb Z$  and it  can be represented as a composition of a radical and an analytic single-valued function: $y(x)=f( (x-p)^{1/m}))$ where $m$ is the  ramification order of $y$. The following corollary follows from Theorem 16.
\proclaim {Corollary 17 ([2], [3])} 1) If a multivalued  function $y$ in the domain $U$   can be  represented by composition of an algebroidal functions of one variable  and single valued analytic functions  then the monodromy group of any ramified germ of $y$  is solvable.

2) If the monodromy group of a ramification germ of $y$ at  $p$ is solvable then in a small neighborhood of $p$ it can be represented by composition of radicals and single valued analytic functions.
\endproclaim

The {\it local monodromy group} of  an algebroidal function $y$ at a point $p\in U$ is the monodromy group of the equation (2) in an arbitrary small neighborhood of the point $p$. The ramified germs of $y$ at the point $p$  correspond to the orbits of the local monodromy group  actions. This statement can be proven in the same way as Lemma 4 was proved.

\subhead {6. Application to the 13-th Hilbert problem}
\endsubhead
In 1957 A.N.~Kolmogorov and V.I.~Arnold proved the following totally unexpected theorem which gave a negative solution to the 13-th Hilbert problem.

\proclaim {Theorem (Kolmogorov--Arnold)} Any continuous function of $n$ variables can be represented as the
composition of functions of a single variable with the help of addition.
\endproclaim

The 13-th Hilbert problem has the following  algebraic version which still remains open:
{\it Is it possible to represent any algebraic function of $n>1$ variables by algebraic functions of a smaller number of variables   with the help of composition and arithmetic operations?}

An {\it entire algebraic function} $y$ in $\Bbb C^N$ is an algebraic  function defined in $U=\Bbb C^N$ by  an equation (2) whose coefficient $f_i$ are polynomials.  An entire algebraic function could be considered as a continuous algebraic function.

It turns out that in Kolmogorov--Arnold Theorem one can not replace  continuous  functions by  entire algebraic functions.

\proclaim {Theorem 18([2], [3])}  If an  entire algebraic function can be represented as a composition of polynomials and entire algebraic functions of one variable, then its local monodromy group at each point is solvable.
\endproclaim

\demo{Proof} Theorem 18 follows from  from Corollary 17.
\enddemo

\proclaim{Corollary 19} A function $y (a,b)$, defined by equation $ y^5 +ay+b=0,$ cannot be expressed in terms of entire algebraic functions of a single variable by means of composition, addition and multiplication.
\endproclaim

\demo{Proof} Indeed, it is easy to check that the local monodromy group of $y$ at the origin is the unsolvable permutation group $S_5$ (see [2], [3]).

\enddemo

Division is not a continuous operation and it destroys the locality. One cannot add  division to the operations used in Theorem 18. It is easy to see that the function $y(a,b)$ from Corollary 18 can be expressed in terms of  entire algebraic functions of a single variable by means of composition  and arithmetic operations: $y(a,b)=g({b}/\root4\of{a^5})\root4\of a$, where $g(u)$ is defined by equation
$g^5+g+u=0$.

The following particular case of the algebraic version of the 13-th Hilbert problem still remains open.

\proclaim{Problem} Show that there is an algebraic function   of two variables which  cannot be expressed in terms of  algebraic functions of a single variable by means of composition  and arithmetic operations.
\endproclaim

\subhead{7. Acknowledgement}\endsubhead
I would like to thank Michael Singer who invited me to write comments for a new edition of the classical J.F.~Ritt's book ``Integration in finite terms''~[5]. This preprint was written as a part of these comments. I also am grateful to Fedor Kogan who edited my English.
\bigskip

\centerline {REFERENCES}
\medskip
[1] V.B. Alekseev, Abel's Theorem in Problems and Solutions. Based on the lectures of Professor
V.I. Arnold, Kluwer Academic Publishers, 2004.

[2] A.G. Khovanskii, The representability of algebroidal functions by superpositions of analytic functions and of algebroidal functions of one variable. Functional Analysis and its applications, V.~4, N~2, 74--79, 1970; translation in Funct. Anal. Appl. 4 (1970), no.~2, 152--156.

[3] A.G. Khovanskii, On compositions of analytical functions with radicals, UMN, 26:2
(1971), 213--214.

[4]  A. Khovanskii, Topological Galois theory. Solvability and unsolvability of equations in
finite terms. Translated by Valentina Kiritchenko and Vladlen Timorin. Series: Springer
Monographs in Mathematics. Springer Berlin Heidelberg. 2014, XVIII, 305 pp. 6 illus.

[5] J. Ritt, Integration in finite terms. Liouville's theory of elementary methods, N. Y. Columbia Univ. Press. 1948.

\end